\documentclass[11pt,a4paper,reqno]{amsart}

\usepackage{amsfonts}
\usepackage{graphicx,tikz}
\usepackage{graphics}
\usepackage{epsfig}
\usepackage{amssymb}
\usepackage{anysize}
\usepackage{lipsum}
\usepackage{wrapfig}
\marginsize{3cm}{3cm}{3cm}{3cm}
\newtheorem{theorem}{Theorem}[section]

\newtheorem{corollary}{Corollary}[section]

\newtheorem{obs}[theorem]{Observation}

\makeatletter

\begin{document}

\setcounter{page}{1}
\vspace{2cm}
\author[\hspace{0.7cm}\centerline{}]{D . Laavanya$^1$, S . Devi Yamini$^2$ }
\title[\centerline{D. Laavanya , S. Devi Yamini : Graceful Coloring of Ladder Graphs
\hspace{0.5cm}}]{Graceful Coloring of Ladder Graphs}

\thanks{\noindent $^1$$^,$$^2$ Vellore Institute of Technology - Department of Mathematics - Kelambakkam - Vandalur Rd, Rajan Nagar, Chennai, Tamil Nadu - 600127, India\\
\indent \,\,\, e-mail: laavanya.d2020@vitstudent.ac.in\\
\indent \,\,\, e-mail: deviyamini.s@vit.ac.in\\}

\begin{abstract}
A graceful k-coloring of a non-empty graph $G=(V,E)$ is a proper vertex coloring $f:V(G)\rightarrow\lbrace 1,2,...,k \rbrace$, $k\geq 2$, which induces a proper edge coloring $f^{*}:E(G)\rightarrow\lbrace 1, 2, . . . , k-1 \rbrace $ defined by $f^{*}(uv) = |f(u)-f(v)|$, where $u,v\in V(G)$. The minimum $k$ for which $G$ has a graceful $k$-coloring is called graceful chromatic number, $\chi_{g}(G)$. The graceful chromatic number for a few variants of ladder graphs are investigated in this article.


 
\end{abstract}
\maketitle 
\bigskip
\bigskip
%
\section{Introduction}
All the graphs $G=(V,E)$ discussed in this paper are connected, simple and finite. Graph labeling introduced by Alexander Rosa in 1967 \cite{10}, is an assignment of integers to the vertices, edges (or both) of a graph $G$ subject to certain conditions. Graph labeling and its types are extensively studied in the literature \cite{4}. Among the various labelings, $\beta$-labeling is one of the prominent labeling. It is also referred as graceful labeling by Golomb \cite{5}, which was initiated to solve the famous Ringel conjecture \cite{10}. Graceful labeling has an extensive range of applications in network addressing, coding theory, communication networks, X-ray crystallography, dental arch, etc. 

Let $G=(V,E)$ be a graph with $m$ edges. An injective function $f:V(G)\rightarrow \lbrace 0,1,2,...m \rbrace$ is a graceful labeling if it induces a bijective function $f^{*}: E(G)\rightarrow \lbrace 1,2,...,m \rbrace$ with the property that for every edge $xy \in E(G)$, $f^{*}(xy)=|f(x)-f(y)| $. If there exists a graceful labeling for a graph $G$, then $G$ is a graceful graph.

A proper coloring of a graph $G$ is an assignment of colors to the vertices or edges of the graph such that every pair of adjacent vertices or edges receive distinct colors respectively. Chromatic number $(\chi(G))$ is the least number of colors required for proper coloring the vertices of the graph $G$, whereas the chromatic index $(\chi^{'}(G))$ is the least number of colors needed for proper coloring the edges of the graph. In \cite{6}, the existence of graceful graphs with arbitrarily large chromatic number was proved.

As an extension of graceful labeling, the concept of graceful chromatic number was introduced by Gary Chartrand in 2015 \cite{2}.  A graceful $k$-coloring of a non-empty graph $G=(V,E)$ is a proper vertex coloring $f:V(G)\rightarrow\lbrace 1,2,...,k \rbrace$, $k\geq 2$, which induces a proper edge coloring $f^{*}:E(G)\rightarrow\lbrace 1, 2, . . . , k-1 \rbrace $ defined by $f^{*}(uv) = |f(u)-f(v)|$, where $u,v\in V(G)$. The minimum $k$ for which $G$ has a graceful $k$-coloring is called graceful chromatic number, $\chi_{g}(G)$.\\

In the introductory paper \cite{2} on graceful coloring, the graceful chromatic number for some well known graphs were computed.
\begin{theorem}\cite{2} \label{1.1}
For a cycle $C_n$, $n\geq 4$, $$\chi_{g}(C_n)=\begin{cases}
		4,& \text{if } n\neq 5\\
		5,& \text{if } n=5\\
		\end{cases}$$
\end{theorem} 
\begin{theorem}\cite{2}
For a path $P_n$, $n\geq 5$, $\chi_{g}(P_n)=5$.
\end{theorem}
\begin{theorem}\cite{2}
For a wheel graph $W_n$, $n\geq 6$, $\chi_{g}(W_n)=n$.
\end{theorem}
\begin{theorem}\cite{2}
If $T$ is a tree with maximum degree $\Delta $, then $\chi_{g}(T)\leq \lceil\frac{5\Delta}{3}\rceil$.
\end{theorem}
\begin{theorem}\cite{2}
If $G$ is a complete bipartite graph of order $n\geq 3$, then $\chi_{g}(G)=n$.
\end{theorem}
\begin{theorem}\cite{2} \label{1.6}
If $G$ is a r-regular graph, then $\chi_{g}(G)\geq r+2$, where $r\geq 2$.
\end{theorem}
\begin{theorem}\cite{2} \label{1.7}
For a nontrivial connected graph $G$, $\chi_{g}(G)\geq \Delta+1$.
\end{theorem}
\begin{theorem}\cite{2} \label{1.8}
For a subgraph $G^{'}$ of G, $\chi_{g}(G^{'})\leq \chi_{g}(G)$.
\end{theorem}
\begin{theorem}\cite{2} \label{1.9}
Let $f:V(G)\rightarrow \lbrace 1,2,...,k \rbrace$, $k\geq 2$ be a coloring of a nontrivial connected graph $G$. Then $f$ is a graceful coloring of $G$ if and only if 
\begin{enumerate}
\item[(i)] for each vertex $v$ of $G$, the vertices in the closed neighborhood $N[v]$ of $v$ are assigned distinct colors by $f$ and
\item[(ii)] for each path $(x,y,z)$ of order $3$ in $G$, $f(y)\neq \frac{f(x)+f(z)}{f(y)}$.
\end{enumerate}
\end{theorem}

Let $T_{\Delta ,h}$ denote the rooted tree (root $v$) with every vertex at a distance less than the height $h$ from $v$ having degree $\Delta $ and the remaining vertices are at a distance $h$ from $v$ as leaves \cite{3}.
\begin{theorem}\cite{3}
For each integer $ \Delta \geq 2$, $\chi_{g}(T_{\Delta ,2}) = \lceil \frac{1}{2}(3 \Delta + 1)\rceil$.
\end{theorem}
\begin{theorem}\cite{3}
For each integer $ \Delta \geq 2$, $\chi_{g}(T_{\Delta ,3}) = \lceil \frac{1}{8}(13 \Delta + 1)\rceil$.
\end{theorem}
\begin{theorem}\cite{3}
For each integer $ \Delta \geq 2$, $\chi_{g}(T_{\Delta ,4}) = \lceil \frac{1}{32}(53 \Delta + 1)\rceil$.
\end{theorem}
\begin{theorem}\cite{3}
For $ \Delta \geq 2$, $h \geq 2+ \lfloor \frac{1}{3} \Delta \rfloor$, $\chi_{g}(T_{\Delta ,h}) = \lceil \frac{5}{3} \Delta \rceil$.
\end{theorem}
The graceful chromatic number of caterpillars were investigated along with a characterization in \cite{13}. The graceful chromatic number for some subclasses of the following graphs have been established in the literature: unicyclic graphs\cite{1}; graphs with diameter at least $2$ \cite{7}; regular and irregular graphs \cite{8}.

\section{Preliminaries}
Denote $[a, b]$ as $\lbrace a, a + 1, . . . , b \rbrace$ and $[a]$ as $[1, a]$, where $a,b \in \mathbb{Z}^{+}$ such that $a<b$.
A closed ladder $L_n$, $n\geq 2$ is a graph obtained from two paths $P_n$ with $V(L_n) = \lbrace x_i,y_i:1\leq i\leq n\rbrace $ and $E(L_n) = \lbrace x_{i}x_{i+1},y_{i}y_{i+1}:1\leq i\leq n-1\rbrace $  $ \cup \lbrace x_{i}y_{i}:1\leq i\leq n\rbrace $.
An open ladder $OL_n$, $n\geq 2$ is a graph formed by removing the edges $x_1y_1$ and $x_ny_n$ from the closed ladder $L_n$.
A slanting ladder $SL_n$, $n\geq 2$ is a graph obtained from two paths $P_n$ with $V(SL_n)=\lbrace x_i,y_i:1\leq i\leq n\rbrace $ and $E(SL_n)=\lbrace x_{i}x_{i+1},y_{i}y_{i+1},x_{i}y_{i+1}:1\leq i\leq n-1\rbrace $.
A triangular ladder $TL_n$, $n\geq 2$ is a graph obtained from two paths $P_n$ with $V(TL_n)=\lbrace x_i,y_i:1\leq i\leq n\rbrace $ and $E(TL_n)=\lbrace x_{i}x_{i+1},y_{i}y_{i+1},x_{i}y_{i+1}:1\leq i\leq n-1 \rbrace$ $ \cup \lbrace x_{i}y_{i}:1\leq i\leq n\rbrace $.
An open triangular ladder $O(TL_n)$, $n\geq 2$ is a graph obtained by removing the edges $x_{1}y_{1}$ and $x_ny_n$ from the triangular ladder $TL_n$.
A diagonal ladder $DL_n$, $n\geq 2$ is a graph obtained by adding the edges $x_{i+1}y_i$, $1\leq i\leq n-1$ in $TL_n$.
An open diagonal ladder $O(DL_n)$, $n\geq 2$ is a graph formed by removing the edges $x_1y_1$ and $x_ny_n$ from the diagonal ladder $DL_n$.
A circular ladder graph $CL_n$, $n\geq 2$ is a graph obtained by adding the edges $x_1x_n$ and $y_1y_n$ in the closed ladder $L_n$. These variants of ladder graphs \cite{12} are illustrated in Figure  $1$.\\
The cartesian product $G \Box G^{'}$ of two simple connected graphs $G$ and $G^{'}$ is a graph with vertices  $V(G \Box G^{'})=V(G) \times V(G^{'})$ and two vertices $(a,a^{'})$ and $(b,b^{'})$ in $G \Box G^{'}$ are adjacent if the distance between $a$ and $b$ is $0$; and $a^{'}$ and $b^{'}$ is $1$ or the distance between $a$ and $b$ is $1$; and $a^{'}$ and $b^{'}$ is $0$.
The strong product $G \boxtimes G^{'}$ of two connected simple graphs $G$ and $G^{'}$ is a graph with vertices $V(G \boxtimes G^{'})=V(G) \times V(G^{'})$ and two vertices $(a,a^{'})$ and $(b,b^{'})$ in $G \boxtimes G^{'}$ are adjacent if the distance between $a$ and $b$ is $0$; and $a^{'}$ and $b^{'}$ is $1$ or the distance between $a$ and $b$ is $1$; and $a^{'}$ and $b^{'}$ is $0$ or the distance between both $a$ and $b$; and $a^{'}$ and $b^{'}$ is $1$ \cite{11}. Note that, the cartesian product of $P_n$ with $P_2$; and $C_n$ with $P_2$ is equivalent to $L_n$ and $CL_n$ respectively. Also, the strong product of $P_n$ with $P_2$ results in $DL_n$.

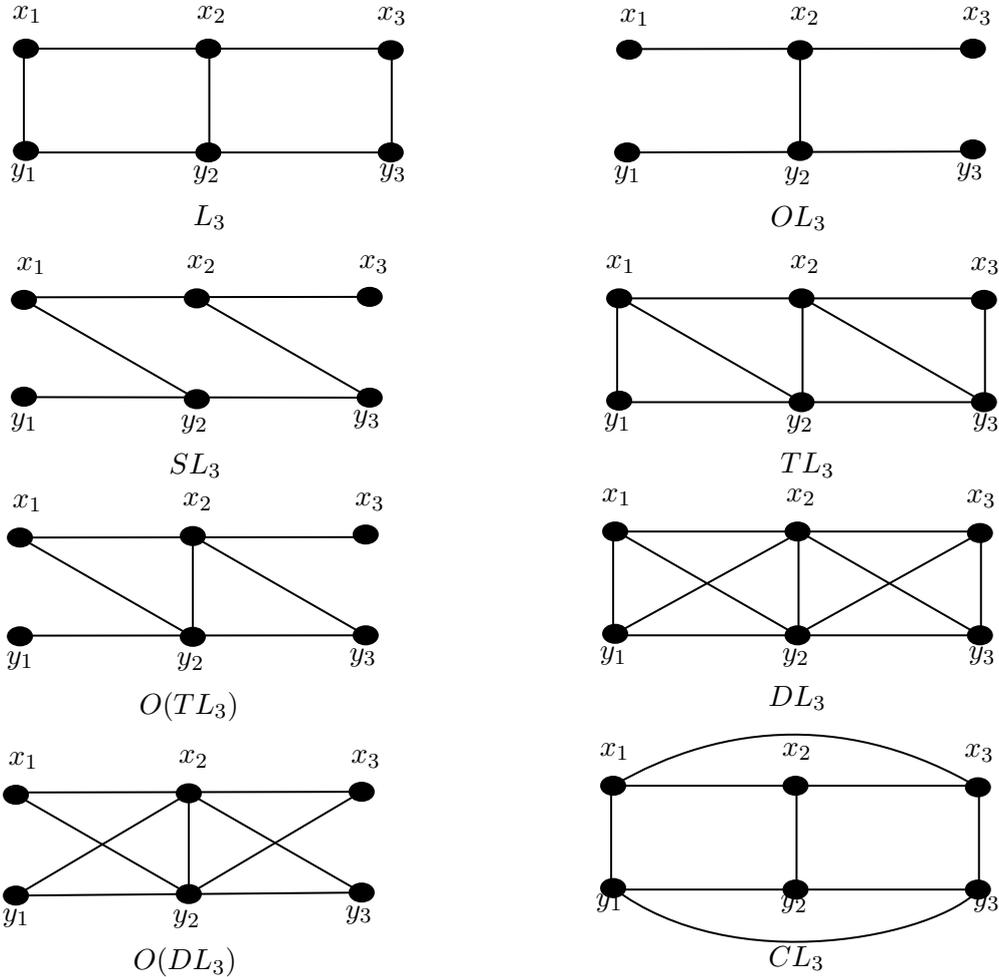
\begin{figure}[ht!]
	\centering
	\tikzset{every picture/.style={line width=0.75pt}} 

\begin{tikzpicture}[x=0.75pt,y=0.75pt,yscale=-1,xscale=1]

\draw   (102,36.48) -- (285.5,36.48) -- (285.5,88.7) -- (102,88.7) -- cycle ;
\draw    (194.5,37.24) -- (194.5,87.94) ;
\draw  [fill={rgb, 255:red, 0; green, 0; blue, 0 }  ,fill opacity=1 ] (97,36.48) .. controls (97,33.97) and (99.69,31.94) .. (103,31.94) .. controls (106.31,31.94) and (109,33.97) .. (109,36.48) .. controls (109,38.99) and (106.31,41.02) .. (103,41.02) .. controls (99.69,41.02) and (97,38.99) .. (97,36.48) -- cycle ;
\draw  [fill={rgb, 255:red, 0; green, 0; blue, 0 }  ,fill opacity=1 ] (188,36.48) .. controls (188,33.97) and (190.69,31.94) .. (194,31.94) .. controls (197.31,31.94) and (200,33.97) .. (200,36.48) .. controls (200,38.99) and (197.31,41.02) .. (194,41.02) .. controls (190.69,41.02) and (188,38.99) .. (188,36.48) -- cycle ;
\draw  [fill={rgb, 255:red, 0; green, 0; blue, 0 }  ,fill opacity=1 ] (97,87.94) .. controls (97,85.43) and (99.69,83.4) .. (103,83.4) .. controls (106.31,83.4) and (109,85.43) .. (109,87.94) .. controls (109,90.45) and (106.31,92.48) .. (103,92.48) .. controls (99.69,92.48) and (97,90.45) .. (97,87.94) -- cycle ;
\draw  [fill={rgb, 255:red, 0; green, 0; blue, 0 }  ,fill opacity=1 ] (188,88.7) .. controls (188,86.19) and (190.69,84.16) .. (194,84.16) .. controls (197.31,84.16) and (200,86.19) .. (200,88.7) .. controls (200,91.2) and (197.31,93.24) .. (194,93.24) .. controls (190.69,93.24) and (188,91.2) .. (188,88.7) -- cycle ;
\draw  [fill={rgb, 255:red, 0; green, 0; blue, 0 }  ,fill opacity=1 ] (279,37.24) .. controls (279,34.73) and (281.69,32.7) .. (285,32.7) .. controls (288.31,32.7) and (291,34.73) .. (291,37.24) .. controls (291,39.74) and (288.31,41.78) .. (285,41.78) .. controls (281.69,41.78) and (279,39.74) .. (279,37.24) -- cycle ;
\draw  [fill={rgb, 255:red, 0; green, 0; blue, 0 }  ,fill opacity=1 ] (279,88.7) .. controls (279,86.19) and (281.69,84.16) .. (285,84.16) .. controls (288.31,84.16) and (291,86.19) .. (291,88.7) .. controls (291,91.2) and (288.31,93.24) .. (285,93.24) .. controls (281.69,93.24) and (279,91.2) .. (279,88.7) -- cycle ;
\draw    (403.89,36.67) -- (575.5,36.48) ;
\draw    (403,88.7) -- (576.33,88.22) ;
\draw    (489.25,37.24) -- (489.25,87.94) ;
\draw  [fill={rgb, 255:red, 0; green, 0; blue, 0 }  ,fill opacity=1 ] (398,36.99) .. controls (398,34.49) and (400.69,32.45) .. (404,32.45) .. controls (407.31,32.45) and (410,34.49) .. (410,36.99) .. controls (410,39.5) and (407.31,41.53) .. (404,41.53) .. controls (400.69,41.53) and (398,39.5) .. (398,36.99) -- cycle ;
\draw  [fill={rgb, 255:red, 0; green, 0; blue, 0 }  ,fill opacity=1 ] (397,88.7) .. controls (397,86.19) and (399.69,84.16) .. (403,84.16) .. controls (406.31,84.16) and (409,86.19) .. (409,88.7) .. controls (409,91.2) and (406.31,93.24) .. (403,93.24) .. controls (399.69,93.24) and (397,91.2) .. (397,88.7) -- cycle ;
\draw  [fill={rgb, 255:red, 0; green, 0; blue, 0 }  ,fill opacity=1 ] (483.25,37.24) .. controls (483.25,34.73) and (485.94,32.7) .. (489.25,32.7) .. controls (492.56,32.7) and (495.25,34.73) .. (495.25,37.24) .. controls (495.25,39.74) and (492.56,41.78) .. (489.25,41.78) .. controls (485.94,41.78) and (483.25,39.74) .. (483.25,37.24) -- cycle ;
\draw  [fill={rgb, 255:red, 0; green, 0; blue, 0 }  ,fill opacity=1 ] (483.25,87.94) .. controls (483.25,85.43) and (485.94,83.4) .. (489.25,83.4) .. controls (492.56,83.4) and (495.25,85.43) .. (495.25,87.94) .. controls (495.25,90.45) and (492.56,92.48) .. (489.25,92.48) .. controls (485.94,92.48) and (483.25,90.45) .. (483.25,87.94) -- cycle ;
\draw  [fill={rgb, 255:red, 0; green, 0; blue, 0 }  ,fill opacity=1 ] (569.5,36.48) .. controls (569.5,33.97) and (572.19,31.94) .. (575.5,31.94) .. controls (578.81,31.94) and (581.5,33.97) .. (581.5,36.48) .. controls (581.5,38.99) and (578.81,41.02) .. (575.5,41.02) .. controls (572.19,41.02) and (569.5,38.99) .. (569.5,36.48) -- cycle ;
\draw  [fill={rgb, 255:red, 0; green, 0; blue, 0 }  ,fill opacity=1 ] (569.5,87.18) .. controls (569.5,84.68) and (572.19,82.64) .. (575.5,82.64) .. controls (578.81,82.64) and (581.5,84.68) .. (581.5,87.18) .. controls (581.5,89.69) and (578.81,91.72) .. (575.5,91.72) .. controls (572.19,91.72) and (569.5,89.69) .. (569.5,87.18) -- cycle ;
\draw    (101.22,161.56) -- (274.5,161.35) ;
\draw    (103,211.78) -- (274.5,212.05) ;
\draw  [fill={rgb, 255:red, 0; green, 0; blue, 0 }  ,fill opacity=1 ] (96,162.86) .. controls (96,160.35) and (98.69,158.32) .. (102,158.32) .. controls (105.31,158.32) and (108,160.35) .. (108,162.86) .. controls (108,165.37) and (105.31,167.4) .. (102,167.4) .. controls (98.69,167.4) and (96,165.37) .. (96,162.86) -- cycle ;
\draw  [fill={rgb, 255:red, 0; green, 0; blue, 0 }  ,fill opacity=1 ] (96,211.56) .. controls (96,209.05) and (98.69,207.02) .. (102,207.02) .. controls (105.31,207.02) and (108,209.05) .. (108,211.56) .. controls (108,214.07) and (105.31,216.1) .. (102,216.1) .. controls (98.69,216.1) and (96,214.07) .. (96,211.56) -- cycle ;
\draw  [fill={rgb, 255:red, 0; green, 0; blue, 0 }  ,fill opacity=1 ] (182.25,162.1) .. controls (182.25,159.59) and (184.94,157.56) .. (188.25,157.56) .. controls (191.56,157.56) and (194.25,159.59) .. (194.25,162.1) .. controls (194.25,164.61) and (191.56,166.64) .. (188.25,166.64) .. controls (184.94,166.64) and (182.25,164.61) .. (182.25,162.1) -- cycle ;
\draw  [fill={rgb, 255:red, 0; green, 0; blue, 0 }  ,fill opacity=1 ] (182.25,212.81) .. controls (182.25,210.3) and (184.94,208.26) .. (188.25,208.26) .. controls (191.56,208.26) and (194.25,210.3) .. (194.25,212.81) .. controls (194.25,215.31) and (191.56,217.35) .. (188.25,217.35) .. controls (184.94,217.35) and (182.25,215.31) .. (182.25,212.81) -- cycle ;
\draw  [fill={rgb, 255:red, 0; green, 0; blue, 0 }  ,fill opacity=1 ] (268.5,161.35) .. controls (268.5,158.84) and (271.19,156.8) .. (274.5,156.8) .. controls (277.81,156.8) and (280.5,158.84) .. (280.5,161.35) .. controls (280.5,163.85) and (277.81,165.89) .. (274.5,165.89) .. controls (271.19,165.89) and (268.5,163.85) .. (268.5,161.35) -- cycle ;
\draw  [fill={rgb, 255:red, 0; green, 0; blue, 0 }  ,fill opacity=1 ] (268.5,212.05) .. controls (268.5,209.54) and (271.19,207.51) .. (274.5,207.51) .. controls (277.81,207.51) and (280.5,209.54) .. (280.5,212.05) .. controls (280.5,214.56) and (277.81,216.59) .. (274.5,216.59) .. controls (271.19,216.59) and (268.5,214.56) .. (268.5,212.05) -- cycle ;
\draw    (102,162.86) -- (188.25,212.81) ;
\draw    (188.25,162.1) -- (274.5,212.05) ;
\draw   (398,162.1) -- (581.5,162.1) -- (581.5,214.32) -- (398,214.32) -- cycle ;
\draw    (490.5,162.86) -- (490.5,213.56) ;
\draw  [fill={rgb, 255:red, 0; green, 0; blue, 0 }  ,fill opacity=1 ] (393,162.1) .. controls (393,159.59) and (395.69,157.56) .. (399,157.56) .. controls (402.31,157.56) and (405,159.59) .. (405,162.1) .. controls (405,164.61) and (402.31,166.64) .. (399,166.64) .. controls (395.69,166.64) and (393,164.61) .. (393,162.1) -- cycle ;
\draw  [fill={rgb, 255:red, 0; green, 0; blue, 0 }  ,fill opacity=1 ] (484,162.1) .. controls (484,159.59) and (486.69,157.56) .. (490,157.56) .. controls (493.31,157.56) and (496,159.59) .. (496,162.1) .. controls (496,164.61) and (493.31,166.64) .. (490,166.64) .. controls (486.69,166.64) and (484,164.61) .. (484,162.1) -- cycle ;
\draw  [fill={rgb, 255:red, 0; green, 0; blue, 0 }  ,fill opacity=1 ] (393,213.56) .. controls (393,211.05) and (395.69,209.02) .. (399,209.02) .. controls (402.31,209.02) and (405,211.05) .. (405,213.56) .. controls (405,216.07) and (402.31,218.1) .. (399,218.1) .. controls (395.69,218.1) and (393,216.07) .. (393,213.56) -- cycle ;
\draw  [fill={rgb, 255:red, 0; green, 0; blue, 0 }  ,fill opacity=1 ] (484,214.32) .. controls (484,211.81) and (486.69,209.78) .. (490,209.78) .. controls (493.31,209.78) and (496,211.81) .. (496,214.32) .. controls (496,216.83) and (493.31,218.86) .. (490,218.86) .. controls (486.69,218.86) and (484,216.83) .. (484,214.32) -- cycle ;
\draw  [fill={rgb, 255:red, 0; green, 0; blue, 0 }  ,fill opacity=1 ] (575,162.86) .. controls (575,160.35) and (577.69,158.32) .. (581,158.32) .. controls (584.31,158.32) and (587,160.35) .. (587,162.86) .. controls (587,165.37) and (584.31,167.4) .. (581,167.4) .. controls (577.69,167.4) and (575,165.37) .. (575,162.86) -- cycle ;
\draw  [fill={rgb, 255:red, 0; green, 0; blue, 0 }  ,fill opacity=1 ] (575,214.32) .. controls (575,211.81) and (577.69,209.78) .. (581,209.78) .. controls (584.31,209.78) and (587,211.81) .. (587,214.32) .. controls (587,216.83) and (584.31,218.86) .. (581,218.86) .. controls (577.69,218.86) and (575,216.83) .. (575,214.32) -- cycle ;
\draw    (403.75,164.37) -- (490,214.32) ;
\draw    (494.75,164.37) -- (581,214.32) ;
\draw    (100,282.43) -- (272.33,282.22) ;
\draw    (100.33,331.78) -- (272.5,331.62) ;
\draw  [fill={rgb, 255:red, 0; green, 0; blue, 0 }  ,fill opacity=1 ] (94,282.43) .. controls (94,279.92) and (96.69,277.89) .. (100,277.89) .. controls (103.31,277.89) and (106,279.92) .. (106,282.43) .. controls (106,284.93) and (103.31,286.97) .. (100,286.97) .. controls (96.69,286.97) and (94,284.93) .. (94,282.43) -- cycle ;
\draw  [fill={rgb, 255:red, 0; green, 0; blue, 0 }  ,fill opacity=1 ] (94,332.13) .. controls (94,329.62) and (96.69,327.59) .. (100,327.59) .. controls (103.31,327.59) and (106,329.62) .. (106,332.13) .. controls (106,334.64) and (103.31,336.67) .. (100,336.67) .. controls (96.69,336.67) and (94,334.64) .. (94,332.13) -- cycle ;
\draw  [fill={rgb, 255:red, 0; green, 0; blue, 0 }  ,fill opacity=1 ] (180.25,281.67) .. controls (180.25,279.16) and (182.94,277.13) .. (186.25,277.13) .. controls (189.56,277.13) and (192.25,279.16) .. (192.25,281.67) .. controls (192.25,284.18) and (189.56,286.21) .. (186.25,286.21) .. controls (182.94,286.21) and (180.25,284.18) .. (180.25,281.67) -- cycle ;
\draw  [fill={rgb, 255:red, 0; green, 0; blue, 0 }  ,fill opacity=1 ] (180.25,332.37) .. controls (180.25,329.87) and (182.94,327.83) .. (186.25,327.83) .. controls (189.56,327.83) and (192.25,329.87) .. (192.25,332.37) .. controls (192.25,334.88) and (189.56,336.91) .. (186.25,336.91) .. controls (182.94,336.91) and (180.25,334.88) .. (180.25,332.37) -- cycle ;
\draw  [fill={rgb, 255:red, 0; green, 0; blue, 0 }  ,fill opacity=1 ] (266.5,280.91) .. controls (266.5,278.41) and (269.19,276.37) .. (272.5,276.37) .. controls (275.81,276.37) and (278.5,278.41) .. (278.5,280.91) .. controls (278.5,283.42) and (275.81,285.45) .. (272.5,285.45) .. controls (269.19,285.45) and (266.5,283.42) .. (266.5,280.91) -- cycle ;
\draw  [fill={rgb, 255:red, 0; green, 0; blue, 0 }  ,fill opacity=1 ] (266.5,331.62) .. controls (266.5,329.11) and (269.19,327.08) .. (272.5,327.08) .. controls (275.81,327.08) and (278.5,329.11) .. (278.5,331.62) .. controls (278.5,334.12) and (275.81,336.16) .. (272.5,336.16) .. controls (269.19,336.16) and (266.5,334.12) .. (266.5,331.62) -- cycle ;
\draw    (100,282.43) -- (186.25,332.37) ;
\draw    (186.25,281.67) -- (272.5,331.62) ;
\draw    (186.25,281.67) -- (186.25,332.37) ;
\draw   (396,279.4) -- (579.5,279.4) -- (579.5,331.62) -- (396,331.62) -- cycle ;
\draw    (488.5,280.16) -- (488.5,330.86) ;
\draw  [fill={rgb, 255:red, 0; green, 0; blue, 0 }  ,fill opacity=1 ] (391,279.4) .. controls (391,276.89) and (393.69,274.86) .. (397,274.86) .. controls (400.31,274.86) and (403,276.89) .. (403,279.4) .. controls (403,281.91) and (400.31,283.94) .. (397,283.94) .. controls (393.69,283.94) and (391,281.91) .. (391,279.4) -- cycle ;
\draw  [fill={rgb, 255:red, 0; green, 0; blue, 0 }  ,fill opacity=1 ] (482,279.4) .. controls (482,276.89) and (484.69,274.86) .. (488,274.86) .. controls (491.31,274.86) and (494,276.89) .. (494,279.4) .. controls (494,281.91) and (491.31,283.94) .. (488,283.94) .. controls (484.69,283.94) and (482,281.91) .. (482,279.4) -- cycle ;
\draw  [fill={rgb, 255:red, 0; green, 0; blue, 0 }  ,fill opacity=1 ] (391,330.86) .. controls (391,328.35) and (393.69,326.32) .. (397,326.32) .. controls (400.31,326.32) and (403,328.35) .. (403,330.86) .. controls (403,333.37) and (400.31,335.4) .. (397,335.4) .. controls (393.69,335.4) and (391,333.37) .. (391,330.86) -- cycle ;
\draw  [fill={rgb, 255:red, 0; green, 0; blue, 0 }  ,fill opacity=1 ] (482,331.62) .. controls (482,329.11) and (484.69,327.08) .. (488,327.08) .. controls (491.31,327.08) and (494,329.11) .. (494,331.62) .. controls (494,334.12) and (491.31,336.16) .. (488,336.16) .. controls (484.69,336.16) and (482,334.12) .. (482,331.62) -- cycle ;
\draw  [fill={rgb, 255:red, 0; green, 0; blue, 0 }  ,fill opacity=1 ] (573,280.16) .. controls (573,277.65) and (575.69,275.62) .. (579,275.62) .. controls (582.31,275.62) and (585,277.65) .. (585,280.16) .. controls (585,282.66) and (582.31,284.7) .. (579,284.7) .. controls (575.69,284.7) and (573,282.66) .. (573,280.16) -- cycle ;
\draw  [fill={rgb, 255:red, 0; green, 0; blue, 0 }  ,fill opacity=1 ] (573,331.62) .. controls (573,329.11) and (575.69,327.08) .. (579,327.08) .. controls (582.31,327.08) and (585,329.11) .. (585,331.62) .. controls (585,334.12) and (582.31,336.16) .. (579,336.16) .. controls (575.69,336.16) and (573,334.12) .. (573,331.62) -- cycle ;
\draw    (401.75,281.67) -- (488,331.62) ;
\draw    (492.75,281.67) -- (579,331.62) ;
\draw    (488.5,280.16) -- (397,330.86) ;
\draw    (579,280.16) -- (488.5,330.86) ;
\draw    (98.11,410.67) -- (270.5,410.32) ;
\draw    (98,462.54) -- (270.5,461.02) ;
\draw  [fill={rgb, 255:red, 0; green, 0; blue, 0 }  ,fill opacity=1 ] (92,411.83) .. controls (92,409.33) and (94.69,407.29) .. (98,407.29) .. controls (101.31,407.29) and (104,409.33) .. (104,411.83) .. controls (104,414.34) and (101.31,416.37) .. (98,416.37) .. controls (94.69,416.37) and (92,414.34) .. (92,411.83) -- cycle ;
\draw  [fill={rgb, 255:red, 0; green, 0; blue, 0 }  ,fill opacity=1 ] (92,462.54) .. controls (92,460.03) and (94.69,458) .. (98,458) .. controls (101.31,458) and (104,460.03) .. (104,462.54) .. controls (104,465.04) and (101.31,467.08) .. (98,467.08) .. controls (94.69,467.08) and (92,465.04) .. (92,462.54) -- cycle ;
\draw  [fill={rgb, 255:red, 0; green, 0; blue, 0 }  ,fill opacity=1 ] (178.25,411.08) .. controls (178.25,408.57) and (180.94,406.54) .. (184.25,406.54) .. controls (187.56,406.54) and (190.25,408.57) .. (190.25,411.08) .. controls (190.25,413.58) and (187.56,415.62) .. (184.25,415.62) .. controls (180.94,415.62) and (178.25,413.58) .. (178.25,411.08) -- cycle ;
\draw  [fill={rgb, 255:red, 0; green, 0; blue, 0 }  ,fill opacity=1 ] (178.25,461.78) .. controls (178.25,459.27) and (180.94,457.24) .. (184.25,457.24) .. controls (187.56,457.24) and (190.25,459.27) .. (190.25,461.78) .. controls (190.25,464.29) and (187.56,466.32) .. (184.25,466.32) .. controls (180.94,466.32) and (178.25,464.29) .. (178.25,461.78) -- cycle ;
\draw  [fill={rgb, 255:red, 0; green, 0; blue, 0 }  ,fill opacity=1 ] (264.5,410.32) .. controls (264.5,407.81) and (267.19,405.78) .. (270.5,405.78) .. controls (273.81,405.78) and (276.5,407.81) .. (276.5,410.32) .. controls (276.5,412.83) and (273.81,414.86) .. (270.5,414.86) .. controls (267.19,414.86) and (264.5,412.83) .. (264.5,410.32) -- cycle ;
\draw  [fill={rgb, 255:red, 0; green, 0; blue, 0 }  ,fill opacity=1 ] (264.5,461.02) .. controls (264.5,458.52) and (267.19,456.48) .. (270.5,456.48) .. controls (273.81,456.48) and (276.5,458.52) .. (276.5,461.02) .. controls (276.5,463.53) and (273.81,465.56) .. (270.5,465.56) .. controls (267.19,465.56) and (264.5,463.53) .. (264.5,461.02) -- cycle ;
\draw    (98,411.83) -- (184.25,461.78) ;
\draw    (184.25,411.08) -- (270.5,461.02) ;
\draw    (184.25,411.08) -- (184.25,461.78) ;
\draw    (184.25,411.08) -- (98,462.54) ;
\draw    (270.5,410.32) -- (184.25,461.78) ;
\draw   (395,407.29) -- (578.5,407.29) -- (578.5,459.51) -- (395,459.51) -- cycle ;
\draw    (487.5,408.05) -- (487.5,458.75) ;
\draw  [fill={rgb, 255:red, 0; green, 0; blue, 0 }  ,fill opacity=1 ] (390,407.29) .. controls (390,404.79) and (392.69,402.75) .. (396,402.75) .. controls (399.31,402.75) and (402,404.79) .. (402,407.29) .. controls (402,409.8) and (399.31,411.83) .. (396,411.83) .. controls (392.69,411.83) and (390,409.8) .. (390,407.29) -- cycle ;
\draw  [fill={rgb, 255:red, 0; green, 0; blue, 0 }  ,fill opacity=1 ] (481,407.29) .. controls (481,404.79) and (483.69,402.75) .. (487,402.75) .. controls (490.31,402.75) and (493,404.79) .. (493,407.29) .. controls (493,409.8) and (490.31,411.83) .. (487,411.83) .. controls (483.69,411.83) and (481,409.8) .. (481,407.29) -- cycle ;
\draw  [fill={rgb, 255:red, 0; green, 0; blue, 0 }  ,fill opacity=1 ] (390,458.75) .. controls (390,456.24) and (392.69,454.21) .. (396,454.21) .. controls (399.31,454.21) and (402,456.24) .. (402,458.75) .. controls (402,461.26) and (399.31,463.29) .. (396,463.29) .. controls (392.69,463.29) and (390,461.26) .. (390,458.75) -- cycle ;
\draw  [fill={rgb, 255:red, 0; green, 0; blue, 0 }  ,fill opacity=1 ] (481,459.51) .. controls (481,457) and (483.69,454.97) .. (487,454.97) .. controls (490.31,454.97) and (493,457) .. (493,459.51) .. controls (493,462.02) and (490.31,464.05) .. (487,464.05) .. controls (483.69,464.05) and (481,462.02) .. (481,459.51) -- cycle ;
\draw  [fill={rgb, 255:red, 0; green, 0; blue, 0 }  ,fill opacity=1 ] (572,408.05) .. controls (572,405.54) and (574.69,403.51) .. (578,403.51) .. controls (581.31,403.51) and (584,405.54) .. (584,408.05) .. controls (584,410.56) and (581.31,412.59) .. (578,412.59) .. controls (574.69,412.59) and (572,410.56) .. (572,408.05) -- cycle ;
\draw  [fill={rgb, 255:red, 0; green, 0; blue, 0 }  ,fill opacity=1 ] (572,459.51) .. controls (572,457) and (574.69,454.97) .. (578,454.97) .. controls (581.31,454.97) and (584,457) .. (584,459.51) .. controls (584,462.02) and (581.31,464.05) .. (578,464.05) .. controls (574.69,464.05) and (572,462.02) .. (572,459.51) -- cycle ;
\draw    (396,407.29) .. controls (452.5,374) and (516.5,371.72) .. (578,408.05) ;
\draw    (396,458.75) .. controls (456.5,504.91) and (557.5,482.97) .. (578.5,459.51) ;

\draw (94.85,13.5) node [anchor=north west][inner sep=0.75pt]    {$x_{1}$};
\draw (186.85,13.5) node [anchor=north west][inner sep=0.75pt]    {$x_{2}$};
\draw (276.85,14.26) node [anchor=north west][inner sep=0.75pt]    {$x_{3}$};
\draw (93.85,92.96) node [anchor=north west][inner sep=0.75pt]    {$y_{1}$};
\draw (184.85,93.72) node [anchor=north west][inner sep=0.75pt]    {$y_{2}$};
\draw (277.85,92.96) node [anchor=north west][inner sep=0.75pt]    {$y_{3}$};
\draw (185,113.64) node [anchor=north west][inner sep=0.75pt]    {$L_{3}$};
\draw (397.85,15.02) node [anchor=north west][inner sep=0.75pt]    {$x_{1}$};
\draw (482.85,14.26) node [anchor=north west][inner sep=0.75pt]    {$x_{2}$};
\draw (568.85,14.26) node [anchor=north west][inner sep=0.75pt]    {$x_{3}$};
\draw (394.85,93.72) node [anchor=north west][inner sep=0.75pt]    {$y_{1}$};
\draw (479.85,94.47) node [anchor=north west][inner sep=0.75pt]    {$y_{2}$};
\draw (565.85,92.2) node [anchor=north west][inner sep=0.75pt]    {$y_{3}$};
\draw (473,114.15) node [anchor=north west][inner sep=0.75pt]    {$OL_{3}$};
\draw (96.85,139.88) node [anchor=north west][inner sep=0.75pt]    {$x_{1}$};
\draw (181.85,139.12) node [anchor=north west][inner sep=0.75pt]    {$x_{2}$};
\draw (267.85,139.12) node [anchor=north west][inner sep=0.75pt]    {$x_{3}$};
\draw (93.85,218.58) node [anchor=north west][inner sep=0.75pt]    {$y_{1}$};
\draw (178.85,219.34) node [anchor=north west][inner sep=0.75pt]    {$y_{2}$};
\draw (264.85,217.07) node [anchor=north west][inner sep=0.75pt]    {$y_{3}$};
\draw (173,238.5) node [anchor=north west][inner sep=0.75pt]    {$SL_{3}$};
\draw (390.85,139.12) node [anchor=north west][inner sep=0.75pt]    {$x_{1}$};
\draw (482.85,139.12) node [anchor=north west][inner sep=0.75pt]    {$x_{2}$};
\draw (572.85,139.88) node [anchor=north west][inner sep=0.75pt]    {$x_{3}$};
\draw (389.85,218.58) node [anchor=north west][inner sep=0.75pt]    {$y_{1}$};
\draw (480.85,219.34) node [anchor=north west][inner sep=0.75pt]    {$y_{2}$};
\draw (573.85,218.58) node [anchor=north west][inner sep=0.75pt]    {$y_{3}$};
\draw (478,238.26) node [anchor=north west][inner sep=0.75pt]    {$TL_{3}$};
\draw (94.85,259.45) node [anchor=north west][inner sep=0.75pt]    {$x_{1}$};
\draw (179.85,258.69) node [anchor=north west][inner sep=0.75pt]    {$x_{2}$};
\draw (265.85,258.69) node [anchor=north west][inner sep=0.75pt]    {$x_{3}$};
\draw (91.85,338.15) node [anchor=north west][inner sep=0.75pt]    {$y_{1}$};
\draw (176.85,338.91) node [anchor=north west][inner sep=0.75pt]    {$y_{2}$};
\draw (262.85,336.64) node [anchor=north west][inner sep=0.75pt]    {$y_{3}$};
\draw (158,358.58) node [anchor=north west][inner sep=0.75pt]    {$O( TL_{3})$};
\draw (388.85,256.42) node [anchor=north west][inner sep=0.75pt]    {$x_{1}$};
\draw (480.85,256.42) node [anchor=north west][inner sep=0.75pt]    {$x_{2}$};
\draw (570.85,257.18) node [anchor=north west][inner sep=0.75pt]    {$x_{3}$};
\draw (387.85,335.88) node [anchor=north west][inner sep=0.75pt]    {$y_{1}$};
\draw (478.85,336.64) node [anchor=north west][inner sep=0.75pt]    {$y_{2}$};
\draw (571.85,335.88) node [anchor=north west][inner sep=0.75pt]    {$y_{3}$};
\draw (472,355.56) node [anchor=north west][inner sep=0.75pt]    {$DL_{3}$};
\draw (92.85,388.85) node [anchor=north west][inner sep=0.75pt]    {$x_{1}$};
\draw (177.85,388.1) node [anchor=north west][inner sep=0.75pt]    {$x_{2}$};
\draw (263.85,388.1) node [anchor=north west][inner sep=0.75pt]    {$x_{3}$};
\draw (89.85,467.56) node [anchor=north west][inner sep=0.75pt]    {$y_{1}$};
\draw (174.85,468.31) node [anchor=north west][inner sep=0.75pt]    {$y_{2}$};
\draw (260.85,466.04) node [anchor=north west][inner sep=0.75pt]    {$y_{3}$};
\draw (155,487.23) node [anchor=north west][inner sep=0.75pt]    {$O( DL_{3})$};
\draw (387.85,384.31) node [anchor=north west][inner sep=0.75pt]    {$x_{1}$};
\draw (569.85,385.07) node [anchor=north west][inner sep=0.75pt]    {$x_{3}$};
\draw (385.85,459.99) node [anchor=north west][inner sep=0.75pt]    {$y_{1}$};
\draw (478,460.75) node [anchor=north west][inner sep=0.75pt]    {$y_{2}$};
\draw (574,459.99) node [anchor=north west][inner sep=0.75pt]    {$y_{3}$};
\draw (472,486.48) node [anchor=north west][inner sep=0.75pt]    {$CL_{3}$};
\draw (478.85,384.31) node [anchor=north west][inner sep=0.75pt]    {$x_{2}$};

\end{tikzpicture}

\end{figure}

\section{Main Results}
\begin{obs}\rm \label{obs1}
If $[\Delta +i]$, $i\in \mathbb{Z}^{+}$ colors are applied in graceful coloring, then the vertex of maximum degree will receive the first and last $i$ colors from $[\Delta +i]$.
\begin{proof}
Let $w$ be a vertex of maximum degree and let $X=\lbrace$ first $i$ colors and last $i$ colors $\rbrace$. We prove $f(w)\in X$. Suppose on the contrary, let $f(w)=a$, where $a \notin X$. Then the $\Delta$ neighbours of $w$ should be distinctly colored from $\lbrace 1,2,...,a-1,a+1,...,\Delta +i\rbrace$. Hence there exist at least two neighbours $p$ and $q$ of $w$ such that $f(p)=a+u$ and $f(q)=a-u$, $u\in [1,\Delta -1]$, a contradiction to the proper edge coloring $(f^{*}(wp)=f^{*}(wq))$. Hence $f(w)\in X$.
\end{proof}
\end{obs}

\begin{theorem}\rm \label{thm1}
$\chi_{g}(L_n)=
\begin{cases}
4,&n=2\\
5,&n \geq 3
\end{cases}$
\begin{proof}
Let $V(L_n)=\lbrace x_i, y_i,1\leq i \leq n \rbrace$ and $E(L_n)=\lbrace x_ix_{i+1}, y_iy_{i+1}, 1\leq i \leq n-1 \rbrace \cup \lbrace x_iy_i, 1\leq i \leq n \rbrace$. Let $x_ix_{i+1}=e^{'}_{i}$, $y_iy_{i+1}=e^{*}_{i}$, $1\leq i\leq n-1$ and $x_iy_{i}=e_{i}$, $1\leq i \leq n$.\\
\textbf{Case 1} ($n=2$):
 Note that $L_2=C_4$ and hence $\chi_{g}(L_2)=4$, by Theorem \ref{1.1}.\\
\textbf{Case 2} ($n=3$): 
Since $L_2$ is a subgraph of $L_3$, $\chi_{g}(L_3)\geq \chi_{g}(L_2)=4$, by the Theorem \ref{1.8}. We now show that $\chi_{g}(L_3)\neq 4$. Suppose that there exist a graceful $4$-coloring of $L_3$. It is clear from the Observation \ref{obs1}, the vertices of maximum degree are colored using the colors $1$ and $4$. Without loss of generality, let $f(x_2)=1$ and $f(y_2)=4$. Then $f(x_1)=3$ and $f(y_1)=2$. Now, the vertices $x_3$ and $y_3$ can be colored using the colors which are at distance at least 3 from them. Thus, $f(x_3)=2$, and hence $f(y_3)=3$ which is a contradiction to the proper edge coloring ($f^{*}(x_3y_3)=1=f^{*}(y_2y_3)$). Hence $\chi_{g}(L_3)\geq 5$. In addition, we prove $\chi_{g}(L_3)\leq 5$.\\
 Define a proper vertex coloring $f:V(L_3)\rightarrow [1,5]$ as 
$f{(v)}=
\begin{cases}
1,& \text{if } v=x_2\\
2,& \text{if } v=x_3,y_1\\
3,& \text{if } v=x_1\\
4,& \text{if } v=y_3\\
5,& \text{if } v=y_2\\
\end{cases}$\\
which induces a proper edge coloring $f^{*}:E(L_3)\rightarrow [1,4]$ as

$f^{*}(e)=
\begin{cases}
1,& \text{if } e=e_1,e^{'}_{2},e^{*}_{2}\\
2,& \text{if } e=e^{'}_{1},e_{3}\\
3,& \text{if } e=e^{*}_{1}\\
4,& \text{if } e=e_{2}\\
\end{cases}$\\
Consequently, $\chi_{g}(L_3)=5$.\\
\textbf{Case 3} ($n>3$): From the Theorem \ref{1.8}, $\chi_{g}(L_n)\geq \chi_{g}(L_3)=5$, for $n>3$. We show that $\chi_{g}(L_n)\leq 5$ by describing a proper vertex coloring $f:V(L_n)\rightarrow [1,5]$ as
$$f{(v)}=\begin{cases}
		1,& \text{if } v=x_{i}: i\equiv 2(\text{mod  }4),y_{j}: j\equiv 0(\text{mod  }4): 1\leq i,j \leq n\\ 
		2,& \text{if } v=x_{i}: i\equiv 3(\text{mod  }4),y_{j}: j\equiv 1(\text{mod  }4): 1\leq i,j \leq n\\
		3,& \text{if } v=x_1\\
		4,& \text{if } v=x_{i}: i\equiv 0(\text{mod  }4),y_{j}: j\equiv 2(\text{mod  }4): 1\leq i,j \leq n\\
		5,& \text{if } v=x_{i}: i\equiv 1(\text{mod  }4),y_{j}: j\equiv 3(\text{mod  }4): 1\leq i,j\leq n \text{ and } i\neq 1 \\
	\end{cases}$$
which induces $f^{*}:E(L_n)\rightarrow [1,4]$ as

$$f^{*}(e)=\begin{cases}
		1,& \text{if } e=\lbrace e_1\rbrace, \lbrace e^{'}_{i}\rbrace, \lbrace e^{*}_{j}\rbrace: i,j\equiv 0(\text{mod  }2): 1\leq i,j \leq n-1\\ 
		2,& \text{if } e=\lbrace e^{'}_1,e^{'}_{i}\rbrace: i\equiv 3(\text{mod  }4),\lbrace e^{*}_{j}\rbrace: j\equiv 1(\text{mod  }4): 1\leq i,j \leq n-1\\
		3, & \text{if } e=\lbrace e_{k}\rbrace: k\equiv 1(\text{mod  }1): 1\leq k \leq n\\
		4, & \text{if } e=\lbrace e^{'}_{i}\rbrace: i\equiv 1(\text{mod  }4),\lbrace e^{*}_{j}\rbrace: j\equiv 3(\text{mod  }4): 1\leq i,j \leq n-1 \text{ and } i\neq 1 \\
	\end{cases}$$
Hence $\chi_{g}(L_n)=5$, for $n>3$.
\end{proof}
\end{theorem}

\begin{corollary}\rm
$\chi_{g}(OL_n)=5,n >3$
\end{corollary}

\begin{theorem}\rm\label{thm2}
$\chi_{g}(SL_n)= 5, n\geq 4$.
\begin{proof}
Let $SL_n$ be the slanting ladder with the vertex set $V(SL_n)=\lbrace x_i, y_i,1\leq i \leq n \rbrace$ and the edge set $E(SL_n)=\lbrace x_ix_{i+1}, y_iy_{i+1}, x_iy_{i+1}, 1\leq i \leq n-1 \rbrace$. Let $x_ix_{i+1}=e^{'}_{i}$, $y_iy_{i+1}=e^{*}_{i}$, $x_iy_{i+1}=e_{i}$. 
Clearly $L_n$ is a subgraph of $SL_n$, $\chi_{g}(SL_n)\geq \chi_{g}(L_n)=5$, by the Theorem \ref{1.8}. 
Define $f:V(SL_n)\rightarrow [1,5]$ as 
$$f{(v)}=\begin{cases}
		1,& \text{if } v=x_i,y_j: i\equiv 1(\text{mod  }4), j\equiv 0(\text{mod  }4) : 1\leq i,j\leq n\\
		2,& \text{if } v=x_i,y_j: i\equiv 3(\text{mod  }4), j\equiv 2(\text{mod  }4) : 1\leq i,j\leq n\\
		4,& \text{if } v=x_i,y_j: i\equiv 2(\text{mod  }4), j\equiv 1(\text{mod  }4) : 1\leq i,j\leq n\\
		5,& \text{if } v=x_i,y_j: i\equiv 0(\text{mod  }4), j\equiv 3(\text{mod  }4) : 1\leq i,j\leq n\\
		\end{cases}$$
which induces $f^{*}:E(SL_n)\rightarrow [1,4]$ as
 $$f^{*}(e)=\begin{cases}
		1,& \text{if } e=e_i: 1\leq i\leq n-1\\ 
		2,& \text{if } e=e^{'}_{j}, e^{*}_{k}: j\equiv 2(\text{mod  }4), k\equiv 1(\text{mod  }4) : 1\leq j,k\leq n-1\\
		3,& \text{if } e=e^{'}_{j}, e^{*}_{k}: j\equiv 1(\text{mod  }2), k\equiv 0(\text{mod  }2) : 1\leq j,k\leq n-1\\
		4,& \text{if } e=e^{'}_{j}, e^{*}_{k}: j\equiv 0(\text{mod  }4), k\equiv 3(\text{mod  }4) : 1\leq j,k\leq n-1\\
	\end{cases}$$
Therefore, $\chi_{g}(SL_n)\leq 5$ implies $\chi_{g}(SL_n)= 5$, for $n \geq 4.$

\end{proof}
\end{theorem}

\begin{theorem}\rm\label{thm3}
$\chi_{g}(TL_n)=
\begin{cases}
6,&n=3,4\\
7,&n\geq 5
\end{cases}$
\begin{proof}
Let $V(TL_n)=\lbrace x_i, y_i,1\leq i \leq n \rbrace$ and $E(TL_n)=\lbrace x_ix_{i+1}, y_iy_{i+1}, x_iy_{i+1}, 1\leq i \leq n-1 \rbrace \cup \lbrace x_iy_i, 1\leq i \leq n \rbrace$. Let $x_ix_{i+1}=e^{'}_{i}$, $y_iy_{i+1}=e^{*}_{i}$, $x_iy_{i+1}=a_{i}$, $1\leq i \leq n-1$ and $x_iy_{i}=e_{i}$, $1\leq i \leq n$.\\
\textbf{Case 1} ($n=3,4$):
Since the maximum degree of $TL_n$ is $4$, we get $\chi_{g}(TL_n)\geq 5$, by the Theorem \ref{1.7}. We claim that, $\chi_{g}(TL_n)\neq 5$. Suppose on the contrary, $\chi_{g}(TL_n)=5$. It is clear that $f(w) \not\in \lbrace 2,3,4 \rbrace $, where $w$ is a vertex of maximum degree, by the Observation \ref{obs1}. For $n=3$, without loss of generality, let $f(x_2)=1$ and $f(y_2)=5$. Obviously, $f(x_1)\notin \lbrace 1,3,5 \rbrace$ and hence $f(x_1)\in \lbrace 2,4 \rbrace$. Without loss of generality, assume $f(x_1)=2$, then the only choice of color for the vertex $y_1$ is $4$. Now $f(y_3)\notin [1,5]$ (by the Theorem \ref{1.9}), which is a contradiction to the assumption that $\chi_{g}(TL_n)=5$. Same argument holds when $f(x_1)=4$. For $n=4$, an induced subgraph of maximum degree vertices of $TL_n$ form a cycle of length $4$ which can be gracefully colored with four distinct colors, by the Theorem \ref{1.1}. But we have only two colors $\lbrace 1,5 \rbrace$, which is a contradiction to the assumption that $\chi_{g}(TL_n)=5$. Hence, at least $6$ colors are needed for graceful coloring of $TL_n$, for $n=3,4$. Thus $\chi_{g}(TL_n)\geq 6$. 
Define $f:V(TL_n)\rightarrow [1,6]$ as 
$$f{(v)}=\begin{cases}
		1,& \text{if } v=y_2\\ 
		2,& \text{if } v=x_2\\
		3,& \text{if } v=y_1,y_4\\
		4,& \text{if } v=x_1,x_4\\
		5,& \text{if } v=y_3\\
		6,& \text{if } v=x_3\\
	\end{cases}$$
which induces $f^{*}:E(TL_n)\rightarrow [1,4]$ as
 $$f^{*}(e)=\begin{cases}
		1,& \text{if } e=e_1, e_{2}, e_{3}, e_{4}\\ 
		2,& \text{if } e=e^{'}_{1}, e^{'}_{3}, e^{*}_{1}, e^{*}_{3}\\
		3,& \text{if } e=a_1,a_2,a_3\\
		4,& \text{if } e=e^{'}_{2}, e^{*}_{2}\\
	\end{cases}$$
Therefore, $\chi_{g}(TL_n)\leq 6$, implies $\chi_{g}(TL_n)= 6$, for $n=3,4.$\\
\textbf{Case 2} ($n\geq 5$): 
Since $TL_4$ is a subgraph of $TL_n$, $\chi_{g}(TL_n)\geq \chi_{g}(TL_4)=6$ (by the Theorem \ref{1.8}). We show $\chi_{g}(TL_n) \neq 6$. Assume the contrary that, $\chi_{g}(TL_n)=6$. It is clear from the Observation \ref{obs1}, $f(w)\not\in \lbrace 3,4 \rbrace $, $w$ is a vertex of maximum degree. Let $H$ be an induced subgraph of maximum degree vertices in $TL_n$. Note that $L_n$, $n\geq 3$ is also a subgraph of $H$ which cannot be gracefully colored using four colors $\lbrace 1,2,5,6 \rbrace $, by the Theorem \ref{thm1}. Hence at least 7 colors are needed for graceful coloring of $TL_n$. Thus $\chi_{g}(TL_n)\geq 7$. We now define a graceful 7-coloring $f$ of $TL_n$. 
Define $f:V(TL_n)\rightarrow [1,7]$ as        
 $$f{(v)}=\begin{cases}
		1,& \text{if } v=y_{j}: j\equiv 0(\text{mod  }3): 1\leq j \leq n\\ 
		2,& \text{if } v=x_{i}: i\equiv 1(\text{mod  }3): 1\leq i \leq n\\
		3,& \text{if } v=x_{i}: i\equiv 2(\text{mod  }3): 1\leq i \leq n\\
		4,& \text{if } v=y_1\\
		5,& \text{if } v=y_{j}: j\equiv 1(\text{mod  }3): 1\leq j \leq n \text{ and } j\neq 1\\
		6,& \text{if } v=x_{i}: i\equiv 0(\text{mod  }3): 1\leq i \leq n\\
		7,& \text{if } v=y_{j}: j\equiv 2(\text{mod  }3): 1\leq j \leq n\\
	\end{cases}$$
which induces $f^{*}:E(TL_n)\rightarrow [1,6]$ as
 $$f^{*}(e)=\begin{cases}
		1,& \text{if } e=\lbrace e^{'}_{i}\rbrace: i\equiv 1(\text{mod  }3), \lbrace a_{l}\rbrace: l\equiv 0(\text{mod  }3): 1\leq i, l\leq n-1\\ 
		2,& \text{if } e=\lbrace e_{1}\rbrace, \lbrace e^{*}_{j} \rbrace, : j\equiv 1(\text{mod  }3), \lbrace a_{l}\rbrace: l\equiv 2(\text{mod  }3): 1\leq l,j \leq n-1 \\
		\text{ } &\text{ and } j\neq 1\\
		3,& \text{if } e=\lbrace e^{*}_{1}\rbrace, \lbrace e^{'}_{i}\rbrace: i\equiv 2(\text{mod  }3),\lbrace e_{k}\rbrace: k\equiv 1(\text{mod  }3): 1\leq i \leq n-1,\\
		\text{ } & 4\leq k\leq n\\
		4,& \text{if } e=\lbrace e^{'}_{i}\rbrace: i\equiv 0(\text{mod  }3), \lbrace e^{*}_{j}\rbrace: j\equiv 0(\text{mod  }3), \lbrace e_{k}\rbrace: k\equiv 2(\text{mod  }3): \\
		\text{  }& 1\leq i, j \leq n-1, 1\leq k\leq n\\
		5,& \text{if } e=\lbrace e_{k}\rbrace: k\equiv 0(\text{mod  }3), \lbrace a_{l}\rbrace: l\equiv 1(\text{mod  }3): 1\leq k \leq n,  1\leq l\leq n-1\\
		6,& \text{if } e=\lbrace e^{*}_{j}\rbrace: j\equiv 2(\text{mod  }3): 1\leq j\leq n-1\\
	\end{cases}$$
Therefore, $\chi_{g}(TL_n)\leq 7$, implies $\chi_{g}(TL_n)= 7$, for $n\geq 5$.
\end{proof}
\end{theorem}
\begin{corollary}\rm
$\chi_{g}(O(TL_n))=7,n\geq 5$
\end{corollary}
\begin{theorem}\rm
$\chi_{g}(DL_n)=
\begin{cases}
8,&n=5,6\\
9,&n \geq 7
\end{cases}$
\begin{proof}
Consider the diagonal ladder $DL_n$ with the vertex and the edge set as follows: $V(DL_n)=\lbrace x_i, y_i,1\leq i \leq n \rbrace$, $E(DL_n)=\lbrace x_ix_{i+1}, y_iy_{i+1}, x_iy_{i+1}, y_ix_{i+1}, 1\leq i \leq n-1 \rbrace \cup \lbrace x_iy_i, 1\leq i \leq n \rbrace $.  Let $x_ix_{i+1}=e^{'}_{i}$, $y_iy_{i+1}=e^{*}_{i}$, $x_iy_{i+1}=a^{'}_{i}$, $y_ix_{i+1}=a^{*}_{i}$, $1\leq i \leq n-1$ and $x_iy_{i}=e_{i}$, $1\leq i \leq n$ \\
\textbf{Case 1} ($n=5,6$):
$TL_n$ is a subgraph of $DL_n$, $\chi_{g}(DL_n) \geq \chi_{g}(TL_n)=7$ (by the Theorem \ref{1.8}). We now show that $\chi_{g}(DL_n)\neq 7$. Suppose on the contrary, $\chi_{g}(DL_n)=7$. Observe that $f(w)\not\in \lbrace 3,4,5 \rbrace $, where $w$ is a vertex of maximum degree (by the Observation \ref{obs1}). Let $H$ be an induced subgraph of maximum degree vertices in $DL_n$. Also $TL_n$ ($n=3,4$) is a subgraph of $H$ which cannot be gracefully colored with $4$ colors $\lbrace 1,2,6,7 \rbrace$ (by the Theorem \ref{thm3}), which implies that our assumption $\chi_{g}(DL_n)=7$, for $n=5,6$ is wrong. Hence at least 8 colors are needed for graceful coloring of $DL_n$. Therefore, $\chi_{g}(DL_n)\geq 8$. \\
Define $f:V(DL_n)\rightarrow [1,8]$ as \\
$f{(v)}=
\begin{cases}
1,& \text{if } v=x_4\\
2,& \text{if } v=y_2,y_5\\
3,& \text{if } v=x_3\\
4,& \text{if } v=y_1,y_6\\
5,& \text{if } v=x_1,x_6\\
6,& \text{if } v=y_4\\
7,& \text{if } v=x_2,x_5\\
8,& \text{if } v=y_3\\
\end{cases}$\\
which induces $f^{*}:E(DL_n)\rightarrow [1,7]$ as
\newline
$f^{*}(e)=
\begin{cases}
1,& \text{if } e=\lbrace e_1,e_6\rbrace, \lbrace a^{'}_{2},a^{'}_{4}\rbrace, \lbrace a^{*}_{2},a^{*}_{4}\rbrace\\
2,& \text{if } e=\lbrace e^{'}_{1},e^{'}_{3},e^{'}_{5}\rbrace, \lbrace e^{*}_{1},e^{*}_{3},e^{*}_{5}\rbrace\\
3,& \text{if } e=\lbrace a^{'}_{1},a^{'}_{3},a^{'}_{5}\rbrace, \lbrace a^{*}_{1},a^{*}_{5}\rbrace\\
4,& \text{if } e=\lbrace e^{'}_{2}\rbrace, \lbrace e^{*}_{4}\rbrace\\
5,& \text{if } e=\lbrace e_{2},e_{3},e_{4},e_{5}\rbrace\\
6,& \text{if } e=\lbrace e^{'}_{4}\rbrace, \lbrace e^{*}_{2}\rbrace\\
7,& \text{if } e=\lbrace a^{*}_{3}\rbrace\\
\end{cases}$\\
Thus, $\chi_{g}(DL_n)\leq 8$. Consequently, $\chi_{g}(DL_n)=8$, for $n=5,6$.\\
\textbf{Case 2} ($n \geq 7$):
Obviously, $\chi_{g}(DL_n)\geq \chi_{g}(DL_6)=8$. We show that, graceful coloring of $DL_n$ need at least 9 colors. Assume the contrary that, $\chi_{g}(DL_n)=8$. It can be seen that $f(w)\not\in \lbrace 4,5 \rbrace $, where $w$ is a vertex of maximum degree (by the Observation \ref{obs1}). Let $H$ be an induced subgraph of maximum degree vertices in $DL_n$. Indeed, $TL_n$, $n\geq 5$ is a subgraph of $H$. By the Theorem \ref{thm3}, the colors $[1,3] \cup [6,8] $ are inadequate for graceful coloring of $DL_n$. Hence, at least 9 colors are required for graceful coloring of $DL_n$. Thus, $\chi_{g}(DL_n)\geq 9$. It remains to show $\chi_{g}(DL_n)\leq 9$ by describing $f:V(DL_n)\rightarrow [1,9]$ as
$$f{(v)}=\begin{cases}
		1,& \text{if } v=x_{i}: i\equiv 2(\text{mod  }4): 1\leq i \leq n\\ 
		2,& \text{if } v=y_{j}: j\equiv 2(\text{mod  }4): 1\leq j \leq n\\
		3,& \text{if } v=y_{j}: j\equiv 0(\text{mod  }4): 1\leq j \leq n\\
		4,& \text{if } v=x_{i}: i\equiv 0(\text{mod  }4): 1\leq i \leq n\\
		5,& \text{if } v=y_{1}\\
		6,& \text{if } v=y_{j}: j\equiv 1(\text{mod  }4): 5\leq j \leq n\\
		7,& \text{if } v=x_{i}: i\equiv 1(\text{mod  }4): 1\leq i \leq n\\
		8,& \text{if } v=x_{i}: i\equiv 3(\text{mod  }4): 1\leq i \leq n\\
		9,& \text{if } v=y_{j}: j\equiv 3(\text{mod  }4): 1\leq j \leq n\\
	\end{cases}$$
which induces $f^{*}:E(DL_n)\rightarrow [1,8]$ as

$$f^{*}(e)=\begin{cases}
		1,& \text{if } e=\lbrace e_{k}\rbrace: k\equiv 1(\text{mod  }1): 1\leq k \leq n\\ 
		2,& \text{if } e=\lbrace e_{1}\rbrace, \lbrace a^{'}_{l}\rbrace: l\equiv 0(\text{mod  }4): 1\leq l \leq n-1\\
		3,& \text{if } e=\lbrace e^{'}_{i}\rbrace: i\equiv 0(\text{mod  }4), \lbrace e^{*}_{1},e^{*}_{j}\rbrace: j\equiv 0(\text{mod  }4) 1\leq i,j \leq n-1\\
		4,& \text{if } e=\lbrace e^{'}_{i}\rbrace :	 i\equiv 3(\text{mod  }4),\lbrace e^{*}_{j}\rbrace: j\equiv 1(\text{mod  }4),\lbrace a^{*}_{1},a^{*}_{m}\rbrace: m\equiv 0(\text{mod  }4):\\
		\text{  }&  1\leq i,j,m \leq n-1 \text{ and } j\neq \lbrace 1 \rbrace\\
		5,& \text{if } e=\lbrace a^{'}_{l}\rbrace: l\equiv 1(\text{mod  }2), \lbrace a^{*}_{m}\rbrace: m\equiv 1(\text{mod  }2):m\neq 1, 1\leq l,m \leq n-1\\
		6,& \text{if } e=\lbrace e^{'}_{i}\rbrace: i\equiv 1(\text{mod  }4), \lbrace e^{*}_{j}\rbrace: j\equiv 3(\text{mod  }4),\lbrace a^{*}_{2},a^{*}_{m}\rbrace: m\equiv 2(\text{mod  }4): \\
		\text{  }&  1\leq i,j,m \leq n-1 \text{ and } m\neq 2\\
		7,& \text{if } e=\lbrace e^{'}_{i}\rbrace: i\equiv 2(\text{mod  }4) , \lbrace e^{*}_{j}\rbrace: j\equiv 2(\text{mod  }4): 1\leq i,j \leq n-1\\
		8,& \text{if } e=\lbrace a^{'}_{l}\rbrace: l\equiv 2(\text{mod  }4): 1\leq l \leq n-1\\
	\end{cases}$$
Hence, $\chi_{g}(DL_n)=9$, for $n\geq 7$.
\end{proof}
\end{theorem}

\begin{corollary}\rm
$\chi_{g}(O(DL_n))=9,n \geq 7$.
\end{corollary}

\begin{theorem}\rm\label{thm2}
For $n\geq 4$,
$\chi_{g}(CL_n)=
\begin{cases}
5,&n\equiv 0(\text{mod  }4)\\
6,&otherwise\\
\end{cases}$

\begin{proof}
A circular ladder $CL_n$ is formed by adding two edges $x_1x_n$ and $y_1y_n$ in the closed ladder $L_n$.\\ \textbf{Case 1} ($n\equiv 0(\text{mod  }4)$): Since $CL_n$ is a $3$-regular graph, $\chi_{g}(CL_n)\geq 5$ (by the Theorem \ref{1.6}). We claim that $\chi_{g}(CL_n)\leq 5$ by defining $f:V(CL_n) \rightarrow [1,5]$ as follows.
For $n=4b+4$, where $b\in \lbrace 0,1,2,3,... \rbrace$  
$$f{(v)}=\begin{cases}
		1,& \text{if } v=\lbrace x_{i}\rbrace: i\equiv 1(\text{mod  }4), \lbrace y_{j}\rbrace: j\equiv 3(\text{mod  }4): 1\leq i,j \leq n\\ 
		2,& \text{if } v=\lbrace x_{i}\rbrace: i\equiv 2(\text{mod  }4), \lbrace y_{j}\rbrace: j\equiv 0(\text{mod  }4): 1\leq i,j \leq n\\
		4,& \text{if } v=\lbrace x_{i}\rbrace: i\equiv 0(\text{mod  }4), \lbrace y_{j}\rbrace: j\equiv 2(\text{mod  }4): 1\leq i,j \leq n\\
		5,& \text{if }v=\lbrace x_{i}\rbrace: i\equiv 3(\text{mod  }4), \lbrace y_{j}\rbrace: j\equiv 1(\text{mod  }4): 1\leq i,j \leq n\\
	\end{cases}$$
which induces $f^{*}:E(CL_n)\rightarrow [1,4]$ 
$$f^{*}(e)=\begin{cases}
		1,& \text{if } e=\lbrace e^{'}_{i}, e^{*}_{j} \rbrace: i,j\equiv 1(\text{mod }2): 1\leq i,j\leq n\\ 
		2,& \text{if } e=\lbrace e_{k} \rbrace: k\equiv 0(\text{mod }2): 1\leq k\leq n\\
		3,& \text{if } e=\lbrace e^{'}_{i}, e^{*}_{j} \rbrace: i,j\equiv 0(\text{mod }2): 1\leq i,j\leq n\\
		4,& \text{if } e=\lbrace e_{k} \rbrace: k\equiv 1(\text{mod }2): 1\leq k\leq n\\
	\end{cases}$$
\textbf{Case 2} ($n\not\equiv 0(\text{mod  }4)$): Obviously, $\chi_{g}(CL_n)\geq \chi_{g}(L_n)=5$, for $n\geq 3$ (by the Theorem \ref{1.8}). It is also clear that, the vertices of $CL_n$ are colored using $\lbrace 1,2,4,5 \rbrace $ (by the Observation \ref{obs1}). We claim that, $\chi_{g}(L_n)\neq 5$. Consider a proper vertex coloring $f$ of $L_n$ as $(1,2,5,4,1,2,5,4,...)$ for the vertices in the upper path and $(5,4,1,2,5,4,1,2,...)$ for the vertices in the lower path; which induces the edge coloring $(1,3,1,3,1,3,1,3,...)$ and  $(1,3,1,3,1,3,1,3,...)$ respectively. Note that, the vertex $x_n$ will not receive the color $4$ in $CL_n$ ($n\not\equiv 0(\text{mod  }4)$).\\
If $f(x_n)=1$, then $f^{*}(x_1x_n)=0$\\
If $f(x_n)=2$, then $f^{*}(x_1x_n)=1=f^{*}(x_1x_2)$\\
If $f(x_n)=5$, then $f^{*}(x_1x_n)=4=f^{*}(x_1y_1)$\\
Note that all the above cases leads to a contradiction to the proper edge coloring. Thus $\chi_{g}(CL_n)\geq 6$. In addition, we prove that $\chi_{g}(CL_n)\leq 6$ by defining $f:V(CL_n) \rightarrow [1,6]$ as follows.
For $n=4b+5$, where $b\in \lbrace 0,1,2,3,... \rbrace$, 
$$f{(v)}=\begin{cases}
		1,& \text{if } v=\lbrace x_{i}\rbrace: i\equiv 1(\text{mod  }4) \text{ and } i\neq n, \lbrace y_{j}\rbrace: j\equiv 3(\text{mod  }4): 1\leq i,j \leq n \\ 
		2,& \text{if } v=\lbrace x_{i}\rbrace: i\equiv 2(\text{mod  }4), \lbrace y_n,y_{j}\rbrace: j\equiv 0(\text{mod  }4) \text{ and } j\neq n-1: \\
		\text{ } &1\leq i,j \leq n \\
		3,& \text{if } v=x_{n-1}\\
		4,& \text{if } v=\lbrace x_n, x_{i}\rbrace: i\equiv 0(\text{mod  }4) \text{ and } i\neq n-1, \lbrace y_{j}\rbrace: j\equiv 2(\text{mod  }4): \\
		\text{ } &1\leq i,j \leq n \\
		5,& \text{if } v=\lbrace x_{i}\rbrace: i\equiv 3(\text{mod  }4), \lbrace y_{j}\rbrace:j\equiv 1(\text{mod  }4)  \text{ and } j\neq n: \\
		\text{ } & 1\leq i,j \leq n \\
		6,& \text{if } v=y_{n-1} \\
		\end{cases}$$
which induces $f^{*}:E(CL_n)\rightarrow [1,5]$ 
 $$f^{*}(e)=\begin{cases}
		1,& \text{if } e=\lbrace e^{'}_{n-1}, e^{'}_{i} \rbrace: i\equiv 1(\text{mod } 2) \text{ and } i\neq n-2, \lbrace e^{*}_{j} \rbrace: j\equiv 1(\text{mod  }2): \\
		\text{ } &1\leq i\leq n, 1\leq j\leq n-4\\ 
		2,& \text{if }  e=\lbrace e^{'}_{n-2}, e_{n} \rbrace, \lbrace e_{k} \rbrace: k\equiv 0(\text{mod  }2) \text{ and } k\neq n-1: 1\leq k \leq n-2 \\
		3,& \text{if }  e=\lbrace e^{'}_{n}, e^{'}_{i} \rbrace: i\equiv 0(\text{mod } 2) \text{ and } i\neq n-1, \lbrace e^{*}_{n},e^{*}_{j} \rbrace: j\equiv 0(\text{mod  }2) \\
		\text{ } &\text{ and } j\neq n-1, e_{n-1}:1\leq i,j\leq n\\ 
		4,& \text{if }  e=e^{*}_{n-1}, \lbrace e_{k} \rbrace: k\equiv 1(\text{mod } 2) \text{ and } k\neq n: 1\leq k\leq n-1\\ 
		5,& \text{if } e= e^{*}_{n-2}\\ 
		\end{cases}$$
For $n=4b+6$, where $b\in \lbrace 0,1,2,3,... \rbrace$, 
$$f{(v)}=\begin{cases}
		1,& \text{if } v=\lbrace x_1,x_4,x_{i}\rbrace: i\equiv 3(\text{mod  }4),\lbrace y_6,y_9,y_{j}\rbrace: j\equiv 1(\text{mod  }4) \text{ and } n\neq 6:\\
		\text{ } &11\leq i \leq n, 13\leq j \leq n\\ 
		2,& \text{if } v=\lbrace x_{2},x_{7} \rbrace, x_{5}: n=6, y_{10} \\
		3,& \text{if } v=x_{5}: n\neq 6, \lbrace x_{i}\rbrace: i\equiv 0(\text{mod  }4),y_{6}: n=6, \lbrace y_3, y_{j}\rbrace: j\equiv 2(\text{mod  }4): \\
		\text{ }&8\leq i\leq n, 14\leq j\leq n\\
		4,& \text{if } v=x_{6}: n=6,\lbrace x_3,x_{i}\rbrace: i\equiv 2(\text{mod  }4),y_{5}: n\neq 6, \lbrace y_{j}\rbrace: j\equiv 0(\text{mod  }4): \\
		\text{ } &10\leq i\leq n, 12\leq j\leq b\\
		5,& \text{if } v=y_{5}:n\neq 6,\lbrace y_{2},y_{7} \rbrace\\
		6,& \text{if } v=\lbrace x_{6},x_{i}\rbrace: i\equiv 1(\text{mod  }4), \lbrace y_{1},y_{4},y_{8},y_{j}\rbrace: j\equiv 3(\text{mod }4): \\
		\text{ } &13\leq i\leq n, 11\leq j\leq n\\
	\end{cases}$$
which induces $f^{*}:E(CL_n)\rightarrow [1,5]$ 
$$f^{*}(e)=\begin{cases}
		1,& \text{if } e=e^{'}_{4}:n=6, \lbrace e^{'}_{1}, e^{'}_{7},e^{'}_{9}\rbrace, e^{*}_{4}:n=6,\lbrace e^{*}_{1},, e^{*}_{7},e^{*}_{9}\rbrace, e_{6}:n=6, \\
		\text{ } &\lbrace e_{3},e_{5},e_{k}\rbrace: k\equiv 0(\text{mod }2): 12\leq k\leq n\\ 
		2,& \text{if } e=e^{'}_{5}:n=6, \lbrace e^{'}_{2}, e^{'}_{4},e^{'}_{8},e^{'}_{i}\rbrace: i\equiv 1(\text{mod }2),e^{*}_{5}:n=6, \lbrace e^{*}_{2}, e^{*}_{4},e^{*}_{j}\rbrace: \\
		\text{ } &j\equiv 1(\text{mod }2), e_{10}:11\leq i,j\leq n\\
		3,& \text{if } e= e^{'}_{6}:n=6,\lbrace e^{'}_{3},e^{'}_{5},e^{'}_{i}\rbrace: i\equiv 0(\text{mod }2), e^{*}_{6}:n=6,\lbrace e^{*}_{3},e^{*}_{5},e^{*}_{j}\rbrace: \\
		\text{  }&j\equiv 0(\text{mod }2), e_{5}:n=6,\lbrace e_{2},e_{7},e_{8} \rbrace: 10\leq i \leq n, 12\leq j\leq n \\
		4,& \text{if } e=e^{'}_{6}:n\neq 6, e^{*}_{6}:n\neq 6, e^{*}_{10},e_{9}:n\neq 6 \\
		5,& \text{if } e= e^{*}_{8}, e_{6}:n\neq 6, \lbrace e_{1},e_{4},e_{k}\rbrace : k\equiv 1(\text{mod }2), 11\leq k\leq n\\
	\end{cases}$$
For $n=4b+7$, where $b\in \lbrace 0,1,2,3,... \rbrace$, 
$$f{(v)}=\begin{cases}
		1,& \text{if } v=\lbrace x_{i}\rbrace: i\equiv 1(\text{mod }4), \lbrace y_{j}\rbrace: j\equiv 3(\text{mod }4):j\neq n: 1\leq i,j \leq n\\ 
		2,& \text{if } v=\lbrace x_{i}\rbrace: i\equiv 2(\text{mod }4) \text{ and } i\neq n-1, \lbrace y_{n}, y_{j}\rbrace: j\equiv 0(\text{mod }4): \\
		\text{ } &1\leq i,j\leq n\\
		3,& \text{if } v=x_{n-1}\\
		4,& \text{if } v=\lbrace x_{n},x_{i}\rbrace: i\equiv 0(\text{mod }4), \lbrace y_{j}\rbrace: j\equiv 2(\text{mod }4) \text{ and } j\neq n-1: \\
		\text{ } &1\leq i,j\leq n\\
		5,& \text{if } v=\lbrace x_{i}\rbrace: i\equiv 3(\text{mod }4) \text{ and } i\neq n, \lbrace  y_{j}\rbrace: j\equiv 1(\text{mod }4):\\
		\text{ } & 1\leq i,j\leq n\\
		6,& \text{if } v=y_{n-1}\\
	\end{cases}$$
which induces $f^{*}:E(CL_n)\rightarrow [1,5]$ 
 $$f^{*}(e)=\begin{cases}
		1,& \text{if } e=\lbrace e^{'}_{n-1}, e^{'}_{i}\rbrace: i\equiv 1(\text{mod }2) \text{ and }i\neq n-2, \lbrace  e^{*}_{j} \rbrace: j\equiv 1(\text{mod }2), \\
		\text{ } &1\leq i,j\leq n\\ 
		2,& \text{if } e= e^{'}_{n-2}, e_{n}, \lbrace  e_{k} \rbrace: k\equiv 0(\text{mod }2) \text{ and } k\neq n-1, 1\leq k\leq n\\
		3,& \text{if } e= \lbrace e^{'}_{n}, e^{'}_{i}\rbrace: i\equiv 0(\text{mod }2) \text{ and }i\neq n-1, \lbrace  e^{*}_{j} \rbrace: j\equiv 0(\text{mod }2) \\
		\text{ } &\text{ and }  j\neq n-1:1\leq i,j\leq n\\
		4,& \text{if } e= e^{*}_{n-1},\lbrace e_{k}\rbrace: k\equiv 1(\text{mod }2) \text{ and }k\neq n, 
		1\leq k\leq n\\
	\end{cases}$$	
Hence, $\chi_{g}(CL_n)=6$, for $n\not\equiv 0(\text{mod  }4)$.  			
\end{proof}
\end{theorem}

\end{document}